\documentclass{amsart}

\usepackage {amsfonts}
\usepackage[english]{babel}
\def\g{\gamma}
\def\G{\Gamma}
\def\d{\delta}
\def\a{\alpha}
\def\b{\beta}
\def\p{\varphi}
\def\e{\varepsilon}
\def\l{\lambda}
\def\L{\Lambda}
\def\o{\omega}
\def\O{\Omega}

\def\D{\mathcal D}
\def\S{\mathcal S}

\def\R{{\mathbb R}}
\def\C{{\mathbb C}}
\def\N{{\mathbb N}}
\def\Z{{\mathbb Z}}
\def\Re{\mbox{Re }}
\def\Im{\mbox{Im }}

\DeclareMathOperator{\Res}{Res}
\DeclareMathOperator{\supp}{supp }
\DeclareMathOperator{\dist}{dist}

\newtheorem*{Th}{Theorem}

\begin{document}

\title{Application of Meyer's theorem on quasicrystals to exponential polynomials and Dirichlet series}

\author{Sergii Yu.Favorov}

\address{Sergii Favorov,
\newline\hphantom{iii}  V.N.Karazin Kharkiv National University
\newline\hphantom{iii} Svobody sq., 4, Kharkiv, Ukraine 61022}
\email{sfavorov@gmail.com}

\maketitle {\small
\begin{quote}
\noindent{\bf Abstract.}
A simple necessary and sufficient condition is given for an absolutely convergent Dirichlet series with imaginary exponents and only real zeros to be a finite product of sines. The proof is based on Meyer's theorem on quasicrystals.

\medskip


\end{quote}
}

   \bigskip

Exponential polynomials of the form
\begin{equation}\label{s}
P(z)=\sum_{1\le j\le N} q_j e^{2\pi i\o_j z}, \qquad q_j\in\C,\quad \o_j\in\R,
\end{equation}
with only real zeros are currently being actively studied. One of the reasons for this is the connection of the zeros of such polynomials with Fourier quasicrystals, proven in the papers of A. Olevsky and A. Ulanovsky \cite{OU},\cite{OU1}. Namely, they showed that the corresponding measure
\begin{equation}\label{a}
\mu_A=\sum_n \d_{a_n},
\end{equation}
where $A=\{a_n\}$ is the zero set of $P$, is a Fourier quasicrystal, and vice versa, the support of each Fourier quasicrystal of the form \eqref{a}
is the zero set of an exponential polynomial of the form \eqref{s}. Then this scheme was generalized to Dirichlet series
\begin{equation}\label{S}
Q(z)=\sum_{\o\in\O} q_\o  e^{2\pi i\o z}, \qquad q_\o\in\C,\quad \o\in\R,\quad\sum_{\o\in\O} |q_\o|<\infty,\quad \O\,\,\text{is bounded},
\end{equation}
(\cite{F2}, \cite{F3}). Previously, such exponential polynomials were used to construct a non-trivial example
  of a Fourier quasicrystal of the form \eqref{a} whose support intersects each arithmetic progression at no more than a finite number of points (P.Kurasov, P.Sarnak \cite{KS}).
 Besides, a simple connection was established between these exponential polynomials and Lee-Yang polynomials of many variables (L.Alon, A.Cohen, K.Vinzant \cite{ACV}).

The simplest example of real-rooted polynomials \eqref{s} is any finite product of sines
\begin{equation}\label{sin}
P_0(z)=Ce^{iaz}\prod_{\j=1}^J\sin(\a_j z+\b_j)\qquad C\in\C,\quad a,\,\a_j,\,\b_j\in\R.
\end{equation}
It is  natural to find conditions for an exponential polynomial or a Dirichlet series to have this form. Obviously, this is precisely the case when $\supp\mu$ is a finite union of arithmetic progressions.

There are many papers (e.g.,\cite{C}, \cite{LO}, \cite{F1}) in which conditions of various kinds have been found  for $\supp\mu$ to have this form. But in these papers, the distances between points of the support of either  $\mu$, or  the Fourier transform $\hat\mu$  are assumed to be uniformly separated from zero.  We use Meyer's result, where the support of $\mu$ is any locally finite measure and $\hat\mu$ is any slowly increasing Radon measure with arbitrary support, and obtain the following theorem:

\begin{Th}
Let $Q(z)$ be Dirichlet series \eqref{S} with  only real zeros,  $h_\g$ be coefficients of the Dirichlet series expansion of the function $Q'(z)/Q(z)$ in the half-planes
\begin{equation} \label{Q}
\begin{gathered}
Q'(z)/Q(z)=\sum_{\g\in\G^*} h_\g e^{2\pi i\g z}+h^+_0,\quad \G^*\subset(0,\infty),\quad (\Im z>0),\\
Q'(z)/Q(z)=\sum_{\g\in\G_*} h_\g e^{2\pi i\g z}+h^-_0,\quad \G_*\subset(-\infty,0),\quad (\Im z<0).
\end{gathered}
\end{equation}
Then the condition
\begin{equation}\label{Qr}
\sum_{|\g|<r} |h_\g|=O(r)\quad (r\to\infty),
\end{equation}
is necessary and sufficient for the representation \eqref{sin}.
\end{Th}
{\bf Remark}. Representations \eqref{Q} for Dirichlet series with only real zeros will be shown below.
\medskip

{\bf Necessity}. It is easily seen that for the function $\sin(\a z+\b),\,\a>0,$ and $\Im z>0$ we have
$$
\frac{\cos(\a z+\b)}{\sin(\a z+\b)}=\sum_{n=1}^\infty(-2i)e^{2i\b n}e^{2i\a nz}-i.
$$
There is a similar representation for $\Im z<0$, and the condition \eqref{Qr} is satisfied.

For $P_0(z)$ from \eqref{sin} we obtain
$$
 \frac{P_0'(z)}{P_0(z)}=[\log P_0(z)]'= ia+\sum_{j+1}^J \a_j\frac{\cos(\a_j z+\b_j)}{\sin(\a_j z+\b_j)}.
$$
Therefore, in this case the condition \eqref{Qr} is also satisfied.

\smallskip
{\bf Sufficiency}. Recall that a continuous function $g(z)$ on the open strip
$$
S_{(\a,\b)}=\{z=x+iy\in\C:\,-\infty\le\a<y<\b\le+\infty\}
$$
is almost periodic  if for any $\a',\,\b',\,\a<\a'<\b'<\b$ and
  $\e>0$ the set of $\e$-almost periods
  $$
E_{\a',\,\b',\e}= \{\tau\in\R:\,\sup_{x\in\R,\a'\le y\le\b'}|g(x+\tau+iy)-g(x+iy)|<\e\}
  $$
is relatively dense, i.e., $E_{\a',\,\b',\e}\cap(x,x+L)\neq\emptyset$ for all $x\in\R$ and some $L$ depending on $\e,\a',\b'$.

A finite sum or product of almost periodic in $S_{(\a,\b)}$ functions is also almost periodic in $S_{(\a,\b)}$, and for almost periodic in $S_{(\a,\b)}$ function $g(z)$
the function $1/g(z)$ provided that $|g(z)|\ge c>0$ is almost periodic too; every Dirichlet series \eqref{S} is almost periodic
in the whole plane $S_{(-\infty,+\infty)}$ (see \cite{B}, App.II).

For any holomorphic almost periodic function $f$ in $S_{(\a,\b)}$ with zero set $A$ for any $\e>0$, and $\a',\,\b',\, \a<\a'<\b'<\b$,
 there exists $m(\e,\a',\b')>0$ such that
 $$
 |f(z)|\ge m(\e,\a',\b')\quad\mbox{for}\quad \a'\le\Im z\le\b', \quad z\not\in A(\e):=\{z:\,\dist(z,A)<\e\}
 $$
 (\cite{L}, Ch.6, Lemma 1), and  numbers
 \begin{equation}\label{num}
\#\{z\in S_{(\a',\b')}\,:\,x\le\Re z\le x+1,\,f(z)=0\}
 \end{equation}
 are bounded uniformly in $x\in\R$ (\cite{L}, Ch.6, Lemma 2).

 Note that each zero of $f$ is counted in \eqref{num} as many times as its multiplicity, hence multiplicities of zeros of every holomorphic almost periodic in $S_{(\a,\b)}$ function
 are uniformly bounded in each substrip $S_{(\a',\b')}$.

In our case the distances between the zero set $A$ of $Q$ and every horizontal line $\Im z=y\neq0$ are uniformly separated from zero, hence, $\inf_{x\in\R}|Q(x+iy)|>0$,
and the function $Q'/Q$ is almost periodic in the half-planes $\Im z>0$ and $\Im z<0$.
By \cite{Z}, Ch.6, the function $1/Q(x+iy)$ for every fixed $y\neq0$ expands into an absolutely convergence Fourier series in $x\in\R$,
and so is the function $Q'/Q$, therefore we get representations \eqref{Q} with an absolutely convergent Dirichlet series.   Further, since zeros of $Q$ are in a horizontal strip of bounded width, we get that the numbers
$\o^+:=\sup\G$ and $\o^-:=\inf\G$ belong to $\G$ (\cite{L}, Ch.6, Sec.2). Now it is easy to check that the function $Q'/Q$ tends to  finite limits as $y\to+\infty$ and $y\to-\infty$, hence  $\G^*\subset(0,+\infty)$ and $\G_*\subset(-\infty,0)$ (see, e.g., \cite{Le}, Part II, Ch.1, Sec.5).
\medskip

Denote by $\D$ the space of $C^\infty$-functions with compact support on $\R$.
The Fourier transform is defined for $\p\in\D$ (in general for $\p\in L^1(\R)$) by the equality
$$
   \hat\p(x)=\int_{\R}\p(t)e^{-2\pi i xt}dt.
$$
 The function
\begin{equation}\label{ext}
   \Phi(z)=\int_{\R}\p(t)e^{-2\pi i zt}dt=\widehat{(\p(t)e^{2\pi yt})}(x),\quad z=x+iy,
\end{equation}
 is the entire extension of $\hat\p$. It is easy to check that for every $k\in\N$ and $s>0$ uniformly in $y\in[-s,s]$
\begin{equation}\label{gro}
  \Phi(x+iy)=O(|x|^{-k})\quad(|x|\to\infty).
\end{equation}
The Fourier transform of a measure $\mu$ on $\R$  is defined by setting
\begin{equation}\label{h}
(\hat\mu,\p)=(\mu,\hat\p),\qquad \p\in\D,
\end{equation}
in the general case $\hat\mu$ is a distribution.

The space $\D$ is dense in the Schwartz space $\S$ of all $C^\infty$-functions $\phi$ such that $\sup_{x\in\R}|x|^k\max_{0\le m\le k}|\phi^{(m)}(x)|<\infty$ for all $k$,
 and the Fourier transform maps continuously $\S$ onto $\S$ with respect to the topology on $\S$ (see, e.g.,\cite{R}). If $|\mu|(-r,r)=O(r^k)$ as $r\to\infty$ with some $k<\infty$,
 then $\mu$ is a continuous linear functional on $\S$. Therefore, for such measures \eqref{h} also holds for all $\p\in\S$.

\medskip
Since the numbers \eqref{num} for $f=Q$ are uniformly bounded, each connected component of $A(\e)$  contains
 no segment of length $1$ for sufficiently  small $\e$.
 Hence there are sequences $L_k\to+\infty,\,L'_k\to-\infty$ and $m=m(\e,\eta)>0$ such that
$$
   |Q(x+iy)|>m\quad\text{for}\quad x=L_k \quad\text{or}\quad x=L'_k,\quad |y|\le\eta.
$$
Let $\p\in\D$, and $\Phi(z)$ be defined in \eqref{ext}.
 Consider  integrals of the function $\Phi(z)Q'(z)Q^{-1}(z)$ over the boundaries of  rectangles
$\Pi_k=\{z=x+iy:\,L'_k<x<L_k,\,|y|<\eta\}$.  Taking into account \eqref{gro}, we get that these integrals  tend to the difference
\begin{equation}\label{in}
   \int_{-i\eta-\infty}^{-i\eta+\infty}\Phi(z)Q'(z)Q^{-1}(z)dz-\int_{i\eta-\infty}^{i\eta+\infty}\Phi(z)Q'(z)Q^{-1}(z)dz=:I_*-I^*
\end{equation}
 as $k\to\infty$. On the other hand, they tend to the sum
  \begin{equation}\label{r}
  2\pi i\sum_{\l:Q(\l)=0}\Res_\l \Phi(z)Q'(z)Q^{-1}(z)=2\pi i\sum_{\l:Q(\l)=0}a(\l)\Phi(\l)=2\pi i(\mu_A,\hat\p),
\end{equation}
where $a(\l)$ is the multiplicity of zero of $Q(z)$ at the point $\l$, and $\mu_A=\sum_{\l:Q(\l)=0}\d_\l$, where every zero $\l$ appears $a(\l)$ times.

Using \eqref{Q} and taking into account  \eqref{gro}, we obtain
$$
 I_*=\sum_{\g\in\G_*}h_\g e^{2\pi\g\eta}\int_{-\infty}^{+\infty}\Phi(x-i\eta)e^{2\pi i\g x}dx+h_0^-\int_{-\infty}^{+\infty}\Phi(x-i\eta)dx,
$$
$$
 I^*=\sum_{\g\in\G^*}h_\g e^{-2\pi\g\eta}\int_{-\infty}^{+\infty}\Phi(x+i\eta)e^{2\pi i\g x}dx+h_0^+\int_{-\infty}^{+\infty}\Phi(x+i\eta)dx,
$$
Therefore the formula for the inverse Fourier transform and \eqref{ext} yield
$$
  I_*=\sum_{\g\in\G_*\cap\supp\p}h_\g\p(\g)+h_0^-\p(0), \qquad I^*=\sum_{\g\in\G^*\cap\supp\p}h_\g\p(\g)+h_0^+\p(0).
$$
Taking into account the definition of $\hat\mu_A$, we get from  \eqref{in} and \eqref{r}
$$
(\hat\mu_A,\p)=(\mu_A,\hat\p)=\sum_{\g\in\G^*\cap\supp\p}\frac{ih_\g}{2\pi}\p(\g)-\sum_{\g\in\G_*\cap\supp\p}\frac{ih_\g}{2\pi}\p(\g)+\frac{i(h_0^+-h_0^-)}{2\pi}\p(0).
$$
By \eqref{Qr}, for any $r\in(1,\infty)$ and all $\p\in\D$ with support in $[-r,r]$ we have
$$
  |(\hat\mu_A,\p)|\le Cr\sup_\R |\p(t)|.
$$
The distribution $\hat\mu_A$ has an expansion to a linear functional with the same bound on the space of continuous functions $g$ on $[-r,r]$ such that $g(-r)=g(r)=0$.
Consequently, $\hat\mu_A$ is a measure of the form
 \begin{equation}\label{ha}
\hat\mu_A=\sum_{\g\in\G^*}\frac{ih_\g}{2\pi}\d_\g-\sum_{\g\in\G_*}\frac{ih_\g}{2\pi}\d_\g+\frac{i(h_0^+-h_0^-)}{2\pi}\d_0,
\end{equation}
that satisfies the condition
$$
|\hat\mu_A|(-r,r)=O(r)\qquad(r\to\infty).
$$
Since  multiplicities of zeros of $Q$ are uniformly bounded by a number $K$, we see that masses of the measure $\mu_A$ take a finite number of values.
Now we use the following  Meyer's theorem:

\begin{Th}[ \cite{M}, p.26, and \cite{KL}]
If $\mu$ is a complex measure with locally finite support such that its masses take a finite number of values, and $\hat\mu$ is a measure such that
$$
|\hat\mu|(-r,r)=O(r)\qquad (r\to\infty),
$$
then each set $\L_s=\{\l\in\supp\mu:\,\mu(\l)=s\}$ belongs to the coset ring\footnote{The coset ring of an abelian group $G$ is the smallest collection of subsets of $G$ which is closed under finite unions,
finite intersections, and complements  and which contains all shifts of all subgroups of $G$.} of $\R$.
\end{Th}
Hence each set $\L_k=\{\l:\,a_\l=k\},\,1\le k\le K$, belongs to the coset ring of $\R$. H.Rosenthal \cite{Ro} proved
that every discrete (without finite limit points) set $\L$  belonging to the coset ring is a finite union of arithmetic progressions up to a finite set. So
\begin{equation}\label{p}
  \L_k=\left[\bigcup_{m=1}^{M_k} L_{k,m}\cup F_k^+\right]  \setminus F_k^-,\quad L_{k,m}=\{a_{k,m}n+b_{k,m}:\,n\in\Z\},\quad 1\le k\le K,\,\,1\le m\le M_k,
\end{equation}
where $F_k^+,\, F_k^-$ are finite. If the intersection of two progressions $L=\{an+b:\,n\in\Z\}$ and $L'=\{a'n+b':\,n\in\Z\}$ contains more than one point, then $a/a'$ is rational and all the sets
$L,\, L',\, L\setminus L',\, L'\setminus L$ are finite unions
of disjoint arithmetic progressions of the type $\{aa'n+b_r:\,n\in\Z\}$. Therefore we can transform \eqref{p} such  that  any distinct $L_{k,m}$ and $L_{k'm'}$ will have at most one point in common.
\medskip

Set
$$
\mu_{k,m}=\sum_{\l\in L_{k,m}} k\,\d_\l,\quad 1\le k\le K,\quad 1\le m\le M_k.
$$
For every $\p\in\D$ the convolutions
$$
\int\p(t-\l)\mu_{k,m}(d\l)=\sum_{\l\in\L_{k,m}}\p(t-\l)
$$
 are periodic functions in the variable $t$. By \eqref{ha}, we also have
 \begin{equation}\label{c}
    (\mu_A,\p(t-\cdot))=(\hat\mu_A,\hat\p e^{2\pi it})=\sum_{\g\in\G^*}\frac{ih_\g}{2\pi}\hat\p(\g)e^{2\pi it\g}-\sum_{\g\in\G_*}\frac{ih_\g}{2\pi}\hat\p(\g)e^{2\pi it\g}+\frac{i(h_0^+-h_0^-)}{2\pi}\hat\p(0).
 \end{equation}
It follows from \eqref{Qr} and \eqref{ext} that for some $C<\infty$
$$
R(r):=\sum_{|\g|<r}|h_\g|=O(r),\qquad |\hat\p(x)|<C\max\{1,\,|x|^{-2}\}.
$$
Therefore,
$$
  \sum_{\g\in\G^*\cup\G_*}|h_\g||\hat\p(\g)|<C\left(R(1)+\int_1^\infty r^{-2}R(dr)\right)=2C\int_1^\infty R(r)r^{-3}dr<\infty.
$$
We obtain that \eqref{c} is almost periodic in $t\in\R$. Set
$$
  \nu=\mu_A-\sum_{k=1}^K\sum_{m=1}^{M_k} \mu_{k,m}.
$$
We see that the function $(\nu,\p(t-\cdot))$ is almost periodic as well. On the other hand, its support is bounded because
$$
\supp\nu\subset\bigcup F_k^+\cup F_k^-\cup(L_{k,m}\cap L_{k'm'}),
$$
where the union is taken over all $1\le k,k'\le K,\,1\le m,m'\le M_k$ provided that pairs $(k,m)$ and $(k',m')$ do not coincide.

So $(\nu,\p(t-\cdot))\equiv0$ for all $\p\in\D$, hence, $\nu\equiv0$,
$$
\mu_A=\sum_{k=1}^K\sum_{m=1}^{M_k} \mu_{k,m},
$$
and the function
$$
  D(z)= \frac{Q(z)}{\prod_{k=1}^K\prod_{m=1}^{M_k}\sin^{k}(\pi z/\a_{k,m}-\pi\b_{k,m}/\a_{k,m})}
$$
is entire without zeros. This  function is of exponential type, therefore,  $D(z)=Ce^{iaz}$ (see \cite{L},Ch.1). If $a\not\in\R$, then $D(z)$ is unbounded on the line
$l=\{z=t+i:\,t\in\R\}$. This is impossible, because
$$
\sup_l|Q(z)|<\infty, \qquad \inf_l|\sin(\pi z/\a_{k,m}-\pi\b_{k,m}/\a_{k,m})|>0,\quad k=1,\dots,K,\quad m=1,\dots,M_k.
$$
 Hence we obtain \eqref{sin}.

\end{document}